\documentclass[11pt]{article}
\usepackage{amsmath,amsfonts,geompsfi}

\renewcommand{\bold}[1]{\medskip \noindent {\bf \boldmath #1
                        }\nopagebreak[4]}

\newcommand{\qed}{\nopagebreak[4]\hspace{.2cm} $\square$ \pagebreak[2]\medskip}

\newtheorem{theorem}{Theorem}[section]

\newcommand{\integers}{{\mathbb Z}}

\newcommand{\reals}{{\mathbb R}}

\newcommand{\makefig}[3]{
	\begin{figure}[htbp]
        \refstepcounter{figure}
	\label{#2}
        \begin{center}~
		#3~\\
		\medskip
                {\sf Figure \thefigure.  #1}
        \end{center}
	\medskip
	\end{figure}
}

\renewcommand{\bold}[1]{\medskip \noindent {\bf \boldmath #1
                        }\nopagebreak[4]}

\newcommand{\bdry}{\partial}

\newcommand{\closure}{\overline}
\newcommand{\compos}{\circ}

\newcommand{\disjunion}{\sqcup}

\usepackage{amssymb}
\newcommand{\nullset}{\varnothing}

\newcommand{\st}{\; | \;}         

\newcommand{\zed}{\integers}

\newcommand{\diam}{\mbox{\rm diam}}

\newcommand{\interior}{\mbox{\rm int}}

\newcommand{\Isom}{\mbox{\rm Isom}}

\newcommand{\Mod}{\mbox{\rm Mod}}

\newcommand{\Teich}{\mbox{\rm Teich}}

\newtheorem{prop}[theorem]{Proposition}
\newtheorem{lem}[theorem]{Lemma}

\newtheorem{defn}[theorem]{Definition}

\newtheorem{quest}[theorem]{Question}

\newcommand{\calC}{{\mathcal C}}

\newcommand{\calM}{{\mathcal M}}
\newcommand{\calN}{{\mathcal N}}

\newcommand{\calP}{{\mathcal P}}

\usepackage{euscript}

\newcommand{\eS}{{\EuScript S}}

\DeclareMathOperator\Fix{Fix}

\begin{document}
\title{
{\bf Curvature and rank of Teichm\"{u}ller space}
\vspace{.2in}}
\author{Jeffrey Brock\thanks{Research supported by the NSF.}  \ and
Benson Farb\thanks{Research supported by the NSF and the Sloan Foundation.}
}

\maketitle

\begin{abstract}
Let $S$ be a surface with genus $g$ and $n$ boundary components and
let $d(S) = 3g-3+n$ denote the number of curves in any pants
decomposition of $S$.  We employ metric properties of the graph of
pants decompositions $C_{\bf P}(S)$ prove that the Weil-Petersson
metric on Teichm\"{u}ller space $\Teich(S)$ is Gromov-hyperbolic if
and only if $d(S) \le 2$.  When $d(S) \ge 3$ the Weil-Petersson metric
has {\em higher rank} in the sense of Gromov (it admits a
quasi-isometric embedding of $\reals^k, k\geq 2$); when $d(S) \le 2$
we combine the hyperbolicity of the complex of curves and the relative
hyperbolicity of $C_{\bf P}(S)$ prove Gromov-hyperbolicity.

We prove moreover that $\Teich(S)$ admits no 
geodesically complete Gromov-hyperbolic metric of finite covolume
when $d(S) \ge 3$, and that no complete Riemannian
metric of pinched negative curvature exists on Moduli space $\calM(S)$
when $d(S) \ge 2$.  
\end{abstract}

\section{Introduction}
The Weil-Petersson metric on Teichm\"uller space 
$\Teich(S)$ has many curious
properties.  It is a Riemannian metric with negative
sectional curvature, but its curvatures are not bounded away from zero 
or negative infinity.  It is  geodesically
convex, but it is not complete.  
In this paper we show that in spite of exhibiting negative curvature
behavior,
the Weil-Petersson metric is not coarsely negatively curved except
for topologically simple surfaces $S$.
Our main theorem answers a question of Bowditch
\cite[Question 11.4]{Bestvina:problems}.  
\begin{theorem}
\label{theorem:main}Let $S$ be a compact surface of genus $g$ with $n$
boundary components.  Then the Weil-Petersson metric on $\Teich(S)$ 
is Gromov-hyperbolic if and only if $3g-3+n \le 2$.
\end{theorem}

The constant $d(S)=3g - 3 + n$ is fundamental
in Teichm\"uller theory: it is the complex dimension of the
Teichm\"uller space $\Teich(S)$, or more topologically, the number of
curves in any pair-of-pants decomposition of $S$.
An equivalent formulation of Theorem 
\ref{theorem:main}, then, is that the Weil-Petersson metric is
Gromov-hyperbolic precisely when the interior $\interior(S)$ is a
torus with at most two punctures or a sphere with at most five
punctures. 

Theorem~\ref{theorem:main} exhibits the first example
of a ($\Mod(S)$-invariant) 
metric on a Teichm\"{u}ller space of (real) dimension
greater than $2$ that is Gromov-hyperbolic.  In contrast, these
Teichm\"{u}ller spaces admit no complete Riemannian
metric of pinched negative sectional curvature (Theorem
\ref{theorem:pinched} below).  To summarize, the overlap of the positive
and negative results in this paper give:

\medskip
\noindent {\sl When $d(S) = 2$, the Weil-Petersson metic on $\Teich(S)$ is
Gromov-hyperbolic, yet $\Teich(S)$ admits no
(equivariant) complete, Riemannian metric with pinched negative curvature.}

\medskip

The fact that the Weil-Petersson metric is not Gromov-hyperbolic when
$d(S)\geq 3$ relies heavily on a geometric investigation of a
combinatorial model for the Weil-Petersson metric constructed in
\cite{Brock:wp} (see below).  

\bigskip
\noindent
{\bf Constraints on metrics on \boldmath${\cal M}_{g,n}$. } In
S. Kravetz' 1959 paper \cite{Kravetz:wrong}, it was claimed that the
Teichm\"{u}ller metric on the moduli space ${\cal M}_{g,n}$ of Riemann
surfaces of genus $g\geq 2$ with $n$ punctures is a complete 
metric with negative curvature in the sense of Busemann, i.e.\ that
$\Teich(S_{g,n})$ admits such an equivarient metric.  The proof had an
error, and the result was shown to be false in \cite{Masur:geods}; in
fact Masur-Wolf showed in
\cite{Masur:Wolf:notGH} that the Teichm\"{u}ller metric is not even
Gromov-hyperbolic (see
\cite{McCarthy:Papadopoulous:MW} 
and
\cite{McCarthy:Papadopoulous:visual} 
for two other proofs).  The following theorem shows
that this phenomenon has little to do with the Teichm\"{u}ller
metric; in fact it holds for a broad class of metrics.

We say that a geodesically complete metric space $X$ has {\em finite
volume}  if for each $\epsilon>0$ there 
is no infinite collection of pairwise disjoint $\epsilon$-balls
embedded in $X$.  
\begin{theorem}
\label{theorem:nohyp}
Suppose $d(S) \geq 3$.  Then 
$\Teich(S)$ admits no proper, geodesically complete, 
Gromov-hyperbolic $\Mod(S)$-equivariant path metric with finite covolume.
\end{theorem}
Theorem \ref{theorem:nohyp} applies in particular to the Teichm\"{u}ller
metric, which satisfies the finite volume condition, giving the main
result of \cite{Masur:Wolf:notGH} when $d(S) \ge 3$.
Theorem \ref{theorem:nohyp} also applies to McMullen's K\"ahler
hyperbolic metric $h$ on $\Teich(S)$ constructed in
\cite{McMullen:kahler} as it is complete and Riemannian, and has finite volume
quotient on $\calM_{g,n}$ (it is also quasi-isometric to the
Teichm\"uller metric).

We note that our argument for Theorem~\ref{theorem:nohyp} breaks down
in the case of the Weil-Petersson metric, which is not complete as a
metric space (and in particular not geodesically complete).  Indeed,
with the isometric action of $\Mod(S)$ on the Weil-Petersson metric, Dehn
twists are infinite order {\em elliptic} elements, as they act with
bounded orbits.  The first part of our proof of
Theorem~\ref{theorem:nohyp} is essentially an argument of
McCarthy-Papadopoulos \cite{McCarthy:Papadopoulous:MW}, where the 
theorem is proven under the additional hypothesis that every
pseudo-Anosov element acts as a hyperbolic isometry.  We reproduce the
argument here for completeness and since it is brief; we also extend the 
proof to the punctured case.

\smallskip

In the absence of the finite-volume hypothesis, it is possible to say
something about complete Riemannian metrics of pinched negative
curvature.
\begin{theorem}
\label{theorem:pinched}
If $3g -3 +n \ge 2$ then ${\cal M}_{g,n}$ admits no complete
Riemannian metric of pinched negative sectional curvature.
\end{theorem}

Note that Theorem~\ref{theorem:nohyp} implies 
Theorem~\ref{theorem:pinched} in  the finite volume case for $d(S) \ge
3$ but does not cover the case when $d(S) = 2$, i.e.\ when $\interior(S)$ is
a torus with two punctures or a sphere with five punctures.



\bold{The rank of the Weil-Petersson metric.}  To prove that the
Weil-Petersson metric is not Gromov-hyperbolic when $d(S) \ge 3$ (the
``only if'' part of Theorem~\ref{theorem:main}), we
show it has {\em higher rank} in these cases. 

The {\em rank} of a metric
space $X$ is the maximal dimension $n$ of a {\em quasi-flat} in $X$,
that is a quasi-isometric embedding $\reals^n\to X$.  This notion of
rank was introduced by Gromov (\cite{Gromov:book:asymptotic}, $6.B_2$), and 
agrees with the usual notion of rank of a nonpositively curved 
symmetric space (this follows easily from the quasi-flats theorem 
of \cite{Eskin:Farb} and \cite{Kleiner:Leeb:rigidity}).

There has been a recurring comparison in the literature of
Teichm\"{u}ller space to symmetric spaces of noncompact type,
particularly in terms of the rank one/higher rank dichotomy (see, e.g.,\
\cite{Ivanov:subgroups,Ivanov:mcg2,Farb:Lubotzky:Minsky}).  The
following theorem adds to the list of {\em higher rank} behavior of
Teichm\"{u}ller space.
\begin{theorem}
\label{theorem:rank}
The rank of the Weil-Petersson metric on $\Teich(S)$ is at least 
$d(S)/2$.
\end{theorem}
As Gromov-hyperbolic metric spaces have rank one, Theorem
\ref{theorem:rank} implies one direction of Theorem 
\ref{theorem:main}.  We conjecture that the Weil-Petersson metric on 
$\Teich(S)$ has rank
precisely the integer part of $(d(S)+1)/2$;
the conjecture is supported by the hierarchies machinery of
\cite{Masur:Minsky:CCII}, but bounding the rank from above appears delicate.

\bigskip
\noindent
{\bf Remark. } While Dehn twists about disjoint simple closed curves
produce quasi-flats in the mapping class group (endowed with the word
metric) - see \cite{Farb:Lubotzky:Minsky}, Dehn twist orbits in 
$\Teich(S)$ {\em do not}
generate quasi-flats in either the Weil-Petersson or Teichm\"uller
metrics: indeed, if $\tau \in \Mod(S)$ is a Dehn twist, we have $d_{\rm
WP}(X, \tau^n X) = O(1)$ and $d_{\rm T}(X, \tau^n X) = O(\log(n)).$
Evidently, the appearance of an orbit of a subgroup generated by
commuting Dehn twists in the Teichm\"uller metric is more akin to a
horosphere in a rank-one symmetric space.  (Note that in higher rank
symmetric spaces horospheres are actually quasi-isometrically embedded.)

\bold{Combinatorics of curves on surfaces.}  
The proof of Theorems~\ref{theorem:main} and~\ref{theorem:rank} rely on
important work of Masur and Minsky
\cite{Masur:Minsky:CCI,Masur:Minsky:CCII} on combinatorial complexes
associated to curves on surfaces.

Let $\eS$ denote the set of all isotopy classes of essential,
non-peripheral, simple
closed curves on the surface $S$.  The {\em curve complex} $\calC(S)$ is
the simplicial complex whose vertices are the elements of $\eS$ and whose
$k$-simplices span collections of $k+1$ curves in $\eS$ that can be
realized pairwise disjointly on $S$.  Metrizing each simplex to be the 
standard Euclidean simplex, one obtains a metric on $\calC(S)$.  The main
theorem of \cite{Masur:Minsky:CCI} is that $\calC(S)$ endowed with this
metric is a Gromov-hyperbolic metric space.

Our proof of Theorem \ref{theorem:rank} uses a closely related
combinatorial object.
Consider the graph whose vertices are pair-of-pants decompositions of
$S$, and whose edges join decompositions that differ by a single {\em
elementary move} (see Figure \ref{figure:moves}).
\makefig{Elementary moves on pants
decompositions.}{figure:moves}{\psfig{file=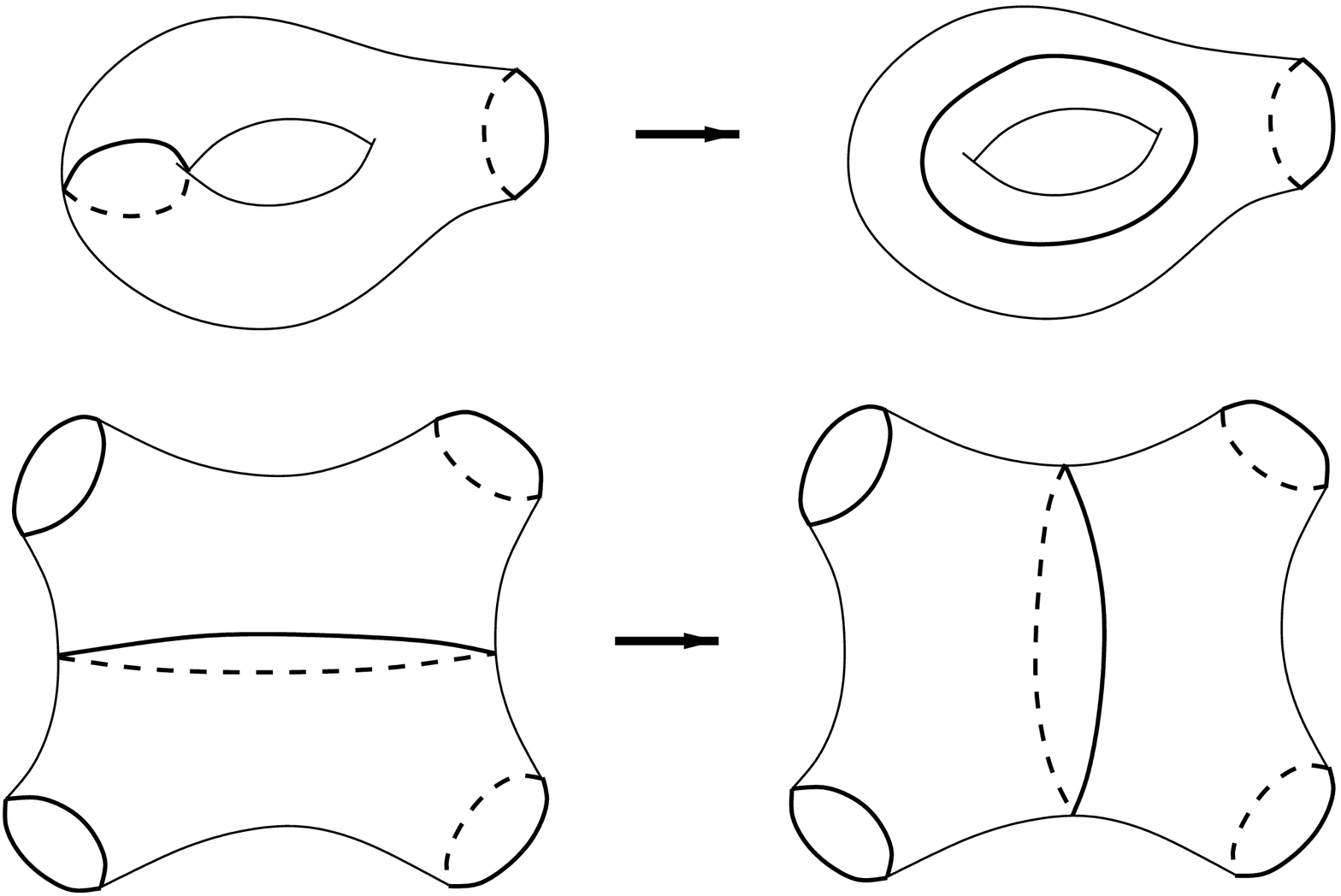,height=2.5in}}
Assigning each edge length 1, we obtain a graph $C_{\bf P}(S)$ with a
distance function $d_{\bf P}(.,.)$ on pairs of vertices given by
taking the minimal length path between two pants-decompositions.  This
graph is the $1$-skeleton of a simplicial complex introduced by
Hatcher-Thurston \cite{Hatcher:Thurston:pants}; in particular they
proved that this graph is connected.  The graph  
$C_{\bf P}(S)$ coarsely models Weil-Petersson geometry.
\begin{theorem}[\cite{Brock:wp}]
The graph $C_{\bf P}(S)$ is quasi-isometric to $\Teich(S)$ endowed with the
Weil-Petersson metric. 
\label{theorem:qi}
\end{theorem}
The proof of Theorem~\ref{theorem:rank} can thus be reduced to finding
quasi-isometrically embedded flats in the graph $C_{\bf P}(S)$.  The
hypothesis is explained by the need for at least two disjoint
essential subsurfaces of $S$ whose Teichm\"uller spaces are themselves
non-trivial, in which case flats arise from pants decompositions
related by elementary moves that occur in disjoint subsurfaces.

Conversely, when $d(S) \le 2$, one cannot perform independent
elementary moves on pairs of pants because there are not enough
disjoint subsurfaces on which to perform them.  The lack of available
subsurfaces leads one to consider the possibility that $C_{\bf P}(S)$
is Gromov-hyperbolicity in these cases.  Gromov-hyperbolicity is known
for $d(S) = 1$ by \cite{Masur:Minsky:CCI}; we establish
Gromov-hyperbolicity of $C_{\bf P}(S)$ for for $d(S) = 2$ using three
ingredients: hyperbolicity of the complex of curves $\calC(S)$ (proved
in \cite{Masur:Minsky:CCI}), the theory of relative hyperbolicity
developed in \cite{Farb:bcp}, and the hierarchical structure of
$\calC(S)$ given in \cite{Masur:Minsky:CCII}.


\bold{Plan of the paper.}
Section \ref{section:prelims} gives preliminaries on Teichm\"uller
theory and Weil-Petersson geometry, Gromov-hyperbolic metric spaces,
and the combinatorics of curves on surfaces. 
Section~\ref{section:nohyp} contains the proof of nonexistence of
geodesically complete Gromov-hyperbolic metrics on $\Teich(S)$ with 
finite covolume when $d(S) \ge 3$ (Theorem~\ref{theorem:nohyp}), and 
gives constraints on complete Riemannian
metrics of pinched negative curvature on ${\cal M}_{g,n}$ 
(Theorem~\ref{theorem:pinched}).  

In section~\ref{section:rank}
we show that the Weil-Petersson metric has higher rank when $d(S) \ge
3$ (Theorem~\ref{theorem:rank}), proving one direction of
Theorem~\ref{theorem:main}.  Finally, in Section~\ref{section:hierarchies}
we prove the other direction of Theorem~\ref{theorem:main} by showing
the Weil-Petersson metric is Gromov-hyperbolic when $d(S) \le 2$.

We conclude the paper with a list of questions for further
investigation.


\bold{Ackhowledgements.}
The authors would like to thank Howard Masur and Yair Minsky for many
informative discussions on the topic of this paper.

\section{Preliminaries}
\label{section:prelims}

\subsection{Teichm\"{u}ller space}

Let $S$ be a topological surface, possibly with boundary.  When $S$
has boundary $\bdry S$ denote by $\interior(S)$ its interior $S -\bdry
S$.  The {\em Teichm\"uller space} $\Teich(S)$ parameterizes finite
area hyperbolic surfaces $X$ equipped with markings, or homeomorphisms
$(f \colon \interior(S) \to X)$ up to isometries that preserve the
marking: i.e.  $$(f \colon \interior(S) \to X) \sim (g \colon
\interior(S) \to Y)$$ if there is an isometry $\phi \colon X \to Y$ so
that $g \simeq \phi
\compos f$.  

The space $\Teich(S)$ is topologized by the distance $$d( (f,X) , (g,
Y)) = \inf_{\varphi} \log L(\varphi)$$ where the infimum is taken over
all bi-Lipschitz diffeomorphisms $\varphi$ isotopic to $g \compos
f^{-1}$ and $L(\varphi)$ is the minimal bi-Lipschitz constant for
$\varphi$.  It is homeomorphic to a cell of dimension $6g -6 +2n$
where $n$ is the number of boundary components of $S$.

The Teichm\"uller space carries a natural complex structure; its
complex cotangent space $T^*_X \Teich(S)$ at $X \in \Teich(S)$ is
identified with the space of {\em holomorphic quadratic differentials} $Q(X)$
on $X$.  The {\em Weil-Petersson metric} on $\Teich(S)$ is obtained by
duality from the $L^2$-inner product on $Q(X)$ $$\langle \varphi ,
\psi \rangle_{\rm WP} = \left( \int_X \frac{\varphi
\closure{\psi}}{\rho^2} |dz|^2 \right)^{1/2}$$
where $\rho(z) |dz|$ is the hyperbolic line element on $X$.

\subsection{The pants complex}
Let $\eS$ denote the isotopy classes of essential simple closed curves
on $S$.  A {\em pants decomposition} $P$ of $S$ is a maximal
collection of distinct elements of $\eS$ so that no two isotopy
classes in $P$ have representatives that intersect.  The {\em pants
complex} $C_{\bf P}(S)$ is the graph with one vertex for each pants
decomposition and an edge joining each pair of vertices whose pants
decompositions differ by an elementary move (see
figure~\ref{figure:moves}).  The distance function $$d_{\bf P} \colon
C^0_{\bf P}(S) \times C^0_{\bf P}(S) \to \zed_{\ge 0}$$ on the vertex
set $C^0_{\bf P}(S)$ of $C_{\bf P}(S)$ counts the minimal number of
elementary moves between maximal partitions.

\bold{Geodesic representatives.}  We briefly describe the
quasi-isometry of Theorem~\ref{theorem:qi}. 

Given $X \in \Teich(S)$, the
elements of $P$ can be represented as pairwise disjoint closed
geodesics on $X$.  The sub-level sets for the lengths of the elements
of $P$ play a special role in Weil-Petersson geometry.  Let
$$V_{L}(P) = \left\{ X \in \Teich(S) \st \ell_X(\alpha) < L \ \
\text{for each} \ \ \alpha \in P  \right\}.$$
Applying theorems of Wolpert \cite{Wolpert:noncompleteness} and Masur \cite{Masur:wp}, 
we have the following (see \cite{Brock:wp}):
\begin{prop}
Given $L$ there is a $D_L$ so that for each $P \in C^0_{\bf P}(S)$, the
sub-level set has Weil-Petersson diameter $$\diam_{\rm WP}(V_L(P)) <D_L.$$
\end{prop}

By a theorem of Bers, there is an $L>0$ so that the union 
$$\bigcup_{P \in C^0_{\bf P}(S)} V_L(P) $$
covers $\Teich(S)$.  Let $V(P) = V_{2L}(P)$.
Let 
$$Q \colon C^0_{\bf P}(S) \to \Teich(S)$$ 
be any map for which $Q(P) \in V(P)$.  Then Theorem~1.1
of \cite{Brock:wp} shows $Q$ is a quasi-isometry.  In other words there are constants
$K_1>1$ and $K_2>0$ depending only on $S$ so that given $X \in V(P_X)$
and $Y \in V(P_Y)$, we have
$$\frac{1}{K_1}d_{\bf P}(P_X,P_Y) - K_2
\le d_{\rm WP}(X,Y) 
\le  d_{\bf P}(P_X, P_Y) + K_2.$$

\subsection{The curve complex}
\label{subsection:cc}

Let $\calC(S)$ be the complex associated to the simple closed curves
$\eS$ on $S$ as follows:
\begin{itemize}
\item The zero-skeleton $\calC^0(S)$ is identified with the elements
of $\eS$.
\item Any $k+1$ curves 
$(\alpha_1,\ldots,\alpha_{k+1})$ in $\eS^{k+1}$ 
with the property $\alpha_i \not= \alpha_j$ and  $i(\alpha_i,
\alpha_j) = 0$, for $i \not= j$ determine a $k$ simplex in $\calC(S)$.
\end{itemize}
For an essential, non-annular subsurface $Y \subset S$, the curve
complex $\calC(Y)$ is a subcomplex of $\calC(S)$.  
The {\em subsurface projection}
$$\pi_Y \colon \calC(S) \to \calP(\calC(Y))$$
from the curve complex $\calC(S)$ to the set $\calP(\calC(Y))$ of all
subsets of $\calC(Y)$ is defined by setting $\pi_Y(\alpha) = \alpha$ if
$\alpha \in \calC(Y)$ and letting 
$$\pi_Y(\alpha) = \cup_{\alpha'}\bdry\calN(\alpha' \cup
\bdry_{\alpha'}Y)$$
be the union over all arcs $\alpha'$ of 
{\em essential} intersection of $\alpha$ with $Y$ of a regular
neighborhood $\bdry\calN(\alpha' \cup
\bdry_{\alpha'}Y)$ of $\alpha' \cup \bdry_{\alpha'}Y$, where
$\bdry_{\alpha'}Y$
represents the components of the boundary of $Y$ that $\alpha'$
intersects (see \cite[Sec. 2]{Masur:Minsky:CCII}).  

For two subsets $A$ and $B$ of $\calC(Y)$, the (semi)-distance $d_Y(A,B)$ is
defined by 
$$d_Y(A,B) = \diam_{\calC(Y)}(A \cup B).$$
This notion of distance allows us to measure the distance between
pants decompositions $P$ and $P'$ {\em relative to $Y$} by letting
$\pi_Y(P) = \cup_{\alpha \in P} \pi_Y(\alpha)$ and letting the {\em
projection distance} $d_Y(P,P')$ between $P$ and $P'$ relative to $Y$
be 
$$d_Y(P,P') = d_Y(\pi_Y(P),\pi_Y(P')).$$

We also record for reference the following Lipschitz property for
the projection $\pi_Y$:
\begin{lem}[Lem. 2.3 of \cite{Masur:Minsky:CCII}]
Let $Y$ be an essential subsurface of $S$.  Then if $\alpha$ and
$\beta$ simple closed curves whose vertices have distance one in
$\calC(S)$ and $\pi_Y(\alpha) \not= \nullset \not= \pi_Y(\beta)$, then
we have $d_Y(\alpha,\beta) \le 2.$
\label{lemma:Lipschitz:simplex}
\end{lem}
For a detailed discussion of the curve complex and related projection
mappings, see \cite{Masur:Minsky:CCI} and \cite{Masur:Minsky:CCII} 

\subsection{Hyperbolic metric spaces}  
\label{subsection:hyperbolic}

In this subsection we briefly recall some material on hyperbolic 
metric spaces.  The standard reference for this material is \cite{GDLH}.  
All statements about hyperbolic metric spaces which we use can be found
in \cite{GDLH}.

A metric space $X$ is {\em proper} if every closed ball in $X$ is
compact.  If for any $x,y \in X$, there exists a parametrized path 
$\gamma:[0,d(x,y)]\rightarrow X$ from $x$ to $y$ with
$d(\gamma(s),\gamma(t))=|s-t|$ for all $s,t\in [0,d(x,y)]$ then 
$X$ is called a {\em geodesic metric space}.  Here we are using the
metric space notion of path length.  

$X$ is a {\em Gromov-hyperbolic}, or {\em $\delta$-hyperbolic} metric
space if there exists $\delta>0$ so that every geodesic triangle $T$ in $X$
is {\em $\delta$-thin}: the $\delta$-neighborhood of any two sides of
$T$ contains the third side.

A $\delta$-hyperbolic metric space has a natural compactification $X\cup 
\partial X$ where $\partial X$ is the collection of Hausdorff
equivalence classes of geodesic rays in $X$.  Every isometry of $X$ acts 
by homeomorphisms on $\partial X$.  We denote the fixed set of the
action of an isometry $g\in \Isom(X)$ in $\partial X$ by $\Fix(g)$.

Combining Theorems~16 and~17 
of Section 8 of \cite{GDLH} gives the following classification theorem.

\begin{theorem}[Classification of isometries]
\label{theorem:classification}
Every isometry $g$ of a $\delta$-hyperbolic metric space is precisely one of 
the following types:
\begin{enumerate}
\item elliptic: every $g$-orbit is bounded.
\item hyperbolic:  $\Fix(g)=\{x,y\}$ for some $x\neq y$ in $\partial X$, 
in which case any $g$-orbit in $X$ is a quasi-geodesic with limit set 
$\{x,y\}$.
\item parabolic: $\Fix(g)=\{x\}$ for some $x\in \partial X$.
\end{enumerate}
\end{theorem}

\subsection{The rank of a metric space} 
A metric space has a {\em quasi-flat of dimension $n$} if there is a
quasi-isometric embedding $F \colon \reals^n \to X$.
We say a metric space has {\em higher rank} if it admits a quasi-flat
of dimension at least 2.  Since $\reals^n$ is not Gromov-hyperbolic,
the quasi-isometric embedding $F$ provides a family of triangles in
$X$ that violates the $\delta$-thin condition for all $\delta >0$.  Thus,
a higher rank metric space is not Gromov-hyperbolic.

\section{Constraints on metrics on ${\cal M}_{g,n}$}
\label{section:nohyp}
In this section we prove Theorems~\ref{theorem:nohyp}
and~\ref{theorem:pinched}.   

The idea of Theorem~\ref{theorem:nohyp} is that a properly discontinuous
action of $\Mod(S)$ on a geodesically complete, $\delta$-hyperbolic
metric space has certain special properties when $S$ is sufficiently
complicated.  Indeed, when $d(S) \ge 3$, then entire group $\Mod(S)$
must act {\em parabolically} with a single parabolic fixed point at
infinity.  This first part is essentially an argument of
McCarthy-Papadopoulos \cite{McCarthy:Papadopoulous:visual}.  As with
isometric actions on hyperbolic space, we show that such an action
cannot have finite volume quotient.

\bold{Notation:} For notational purposes in the following arguments,
we will use $S_{g,n}$ to refer to a surface with genus
$g$ and $n$ boundary components.

\bold{Proof:} {\em (of Theorem \ref{theorem:nohyp}).}
Suppose to the contrary that $X=\Teich(S_{g,n})$ admits a
$\Mod(S_{g,n})$-equivariant, geodesically complete path metric which
is $\delta$-hyperbolic.  We have that $\Mod(S_{g,n})$ acts properly
discontinuously on $X$ by isometries, and (see
\S\ref{subsection:hyperbolic}) that $\Mod(S_{g,n})$
thus acts by homeomorphisms on the Gromov boundary $\partial X$.

We claim that since $3g-3+n\geq 3$, we may pick a generating set $\{g_i\}$ 
for $\Mod(S_{g,n})$ consisting of Dehn 
twists about non-separating curves, with the following properties:
\begin{enumerate}
\item Each $g_i$ is conjugate in $\Mod(S_{g,n})$ to each $g_j$.
\item The group generated by elements commuting with $g_1$ is not
virtually cyclic, i.e.\ it does not
contain a cyclic subgroup of finite index;  similary for $g_2$.
\item $g_1$ and $g_2$ do not commute; in fact sufficiently high powers
of $g_1$ and $g_2$ generate a free group.
\item The {\em commuting graph} for $\{g_i\}$, consisting of a vertex
for each $g_i$ and an edge connecting commuting elements, is connected. 
\end{enumerate}

When $n=0$ one may take the ``standard'' Dehn-Lickorish-Humphries
generators (see \cite{Farb:Franks:homeoI}).  For $n>0$ one proceeds
inductively by  
using the exact sequence (see, e.g.\ \cite{Ivanov:mcg2}): 
$$1\rightarrow \pi_1(S_{g,n})\rightarrow
\Mod(S_{g,n+1})\rightarrow \Mod(S_{g,n})\rightarrow 1$$ 
where the kernel is generated by ``pushing the puncture'' around a given
loop, one for each loop in a standard generating set for
$\pi_1(S_{g,n})$.  Such elements of $\Mod(S_{g,n+1})$ are generated
by elements $\alpha_1\alpha_2^{-1}$ where each $\alpha_i$ is a Dehn
twist about a non-separating curve; one then adds these Dehn twists to
the previous list of generators, easily checking the required
properties, and continues inductively.

The elements $g_1$ and $g_2$ may then be taken to be an intersecting
pair of loops in the Dehn-Humphries-Lickorish generating set.  
In particular, the group generated by elements commuting
with $g_i$ (for $i=1$ or $i=2$), contains the mapping class group 
$\Mod(S_{g,n-1})$, which is not virtually cyclic for $3g-3+n\geq 3$.

Now apply the classification of isometries of $\delta$-hyperbolic metric
spaces (see Theorem~\ref{theorem:classification} above) to 
$g_1$.  Note that this classification uses the geodesic assumption on
the metric space $X$.  As $g_1$ has infinite order and the action of
$\Mod(S_{g,n})$ is properly discontinuous, it follows that $g_1$ is
not of elliptic type.  Suppose $g_1$ is of hyperbolic type.  Then
$\Fix(g_1)=\{x,y\}$ for some $x\neq y$ in $\partial X$.

The subgroup $H$ of $\Mod(S_{g,n})$ commuting with $g_1$ clearly 
leaves $\{x,y\}$ invariant, and is not virtually cyclic.  But 
Theorem 30 in Section 8 of \cite{GDLH} states that, for a group 
$\Gamma$ acting properly discontinuously on a proper,
geodesic, $\delta$-hyperbolic metric space $X$, the stabilizer of a pair of 
distinct points $\{x,y\}$ has a cyclic subgroup of
finite index, a contradiction.  Thus it must be that $g_1$ is of parabolic
type, and so $\Fix(g_1)=x$ for some $x\in \partial X$.

As each $g_i$ is conjugate to $g_1$, each $g_i$ is also of parabolic
type, say fixing the unique point $x_i\in \partial X$.  As Dehn twists
about disjoint curves commute, since $[g,h]=1$ implies
$g(\Fix(h))=\Fix(h)$, and since the commuting graph of $\{g_i\}$ is
connected, it follows that $\Fix(g)=x$ for each $g\in\{g_i\}$.

We claim that every element of $\Mod(S_{g,n})$, not just the generating
set, is of parabolic type
with unique fixed point $x\in \partial X$. As $\Mod(S_{g,n})$ has a
torsion-free subgroup of finite index, if this is not true then there
exists $g\in\Mod(S_{g,n})$ acting on $X$ as an isometry of
hyperbolic type, with $x$ as an attracting point.  Pick any $z\in X$.
By Theorem 8.21 of
\cite{GDLH}, the orbit $\{g^nz: n\in \zed\}$ is a 
$K$-quasigeodesic in $x$ for some $K\geq 1$ and has limit
point $x$ as $n\rightarrow \infty$. Pick any element $h\in \Mod(S_{g,n})$
acting as a parabolic fixing $x$ (such $h$ exist by the previous
paragraph).  Since $\{g^nz: n\in \zed\}$ and $\{hg^nz: n\in \zed\}$ 
are both $K$-quasigeodesics which 
limit to $x$ as $n\rightarrow \infty$, it follows that there exists
$N>0, C>0$ so that $d(g^n(z),hg^n(z))<C$ for all $n>N$.  But then 
$d(g^{-n}hg^n(z),z)=d(hg^n(z),g^n(z))<C$ for all $n>N$.  Since $h$ is of 
parabolic type and $g$ is of hyperbolic type, the set $\{g^{-n}hg^n:
n\in \zed\}$ is infinite, so that the hypothesis of proper discontinuity
is violated.  Hence the claim is proved \footnote{While this claim is 
obviously true in the negatively curved Riemannian context without the proper
discontinuity hypothesis, it may not be true without it for arbitrary 
Gromov-hyperbolic spaces; see \cite{GDLH}, 8.13.}.

It follows that $\Mod(S_{g,n})$ permutes the set of equivalence classes of
geodesic rays with $x$ as their common endpoint at infinity $x\in
\partial X$.  We now show that one of the hypotheses of the theorem must 
be violated.

Recall from \S 8.1 of \cite{GDLH} that any choice of a point $x\in
\partial X$ and basepoint $y\in X$ determines a 
{\em Busemann function} $\beta:X\to \reals$ on $X$ defined by
$$\beta(z)=\sup_\gamma\{\limsup_{t\to\infty}(d_X(z,\gamma(t))-t))\}$$
where the $\sup$ is taken over all geodesic rays $\gamma$ based at $y$
with $\gamma(\infty)=x$.  In the proof of Proposition 8.18 of \cite{GDLH}
(see also Remark 8.13.ii) it is shown that any parabolic
isometry $g$ fixing $x\in \partial X$ must {\em almost preserve} level 
sets of the Busemann function, that is there exists a constant $C$ so that 
$|\beta(w)-\beta(g^n(w))|\leq C$ for all $n\in \zed$ and any $w$ lying on 
the ray $\gamma$.  

As this holds true for all $g\in \Mod(S_{g,n})$, it follows that
any $\Mod(S_{g,n})$-orbit in $X$ lies within a bounded distance of
some level set of $\beta$.  But $\beta$ is clearly proper on $\gamma$;
in particular there exists a constant $\epsilon$ and points $z_i,
i=1,2,\ldots$ on $\gamma$ with the property that for any $i,j$ with 
$i\neq j$, the
$\epsilon$-ball centered at $z_i$ is disjoint from the
$\Mod(S_{g,n})$-orbit of the $\epsilon$-ball centered at $z_j$.  In
particular the quotient $\Teich(S_{g,n})/\Mod(S_{g,n})$
contains an infinite, disjoint collection of $\epsilon$-balls.  This
contradicts the finite volume hypothesis.
\qed

For Theorem~\ref{theorem:pinched}, we illustrate that the
existence of a $\Mod(S)$-equivariant complete Riemannian metric on
$\Teich(S)$ of {\em pinched
negative curvature} puts even stronger restrictions on an isometric
action of $\Mod(S)$.  

\bold{Proof:} {\em (of Theorem \ref{theorem:pinched}).}
If ${\cal M}_{g,n}$ did admit such a metric, then lifting this metric to
the universal cover $\Teich(S_{g,n})$ gives a properly
discontinuous, isometric action of $\Mod(S_{g,n})$ on a complete,
$1$-connected, pinched negatively curved manifold
$X=\Teich(S_{g,n})$. 

Since $3g-3+n\geq 2$, there exists a
subgroup $N$ of $\Mod(S_{g,n})$ generated by Dehn twists as above,
with elements conjugate in $\Mod(S_{g,n})$, and 
with the property that $N$ is not virtually nilpotent (since it
contains, for example, noncyclic free subgroup 
(see, e.g. \cite{Farb:Lubotzky:Minsky})).   
Note that in case $g=0$, the elements $g_i$ can be taken to be Dehn
twists about curves which surround two punctures (hence are separating), 
in particular they are conjugate in $\Mod(S_{g,n})$.

As $X$ is $\delta$-hyperbolic, the exact same argument as in the proof
of Theorem \ref{theorem:nohyp} gives some $x\in X$ which is fixed 
by each generator $g_i$.

Fix any horosphere $H$ based at $x$, and fix a basepoint $s\in H$.  Then
$d(s,g_i s)\leq C_H$ for some constant $C_H$ for each generator $g_i$.
As the sectional curvature of $X$ is pinched between two negative
constants, we may find a horosphere $H$ based at $x\in \partial X$ so
that $C_H$ is as small as we want; the key point is that the pinching of
the curvatures gives a definite (exponential) rate at which geodesic
rays asymptotic to $x$ converge.  Choosing $H$ so that $C_H$ is smaller
than the Margulis constant for $X$ (which depends only on $\dim(X)$ and
on the pinching constants of the sectional curvatures) and applying the
Margulis Lemma (see, e.g.\ \cite{Ballman:Gromov:Schroeder}), it follows that
$N$ contains a nilpotent subgroup of finite index, a contradiction.
\qed

\bigskip
\noindent
{\bf Remark. }Theorem \ref{theorem:pinched} may also be deduced (though
not in some of the low genus cases) from the topological structure of
the end of ${\cal M}_{g,n}$.  On the one hand, every end of a 
finite volume manifold of pinched negative curvature is homeomorphic to
the product of a compact nilmanifold and $[0,\infty)$; this is
essentially the Margulis Lemma (see
\cite{Ballman:Gromov:Schroeder}). 
On the other hand, it seems to be well-known (see, e.g.\
\cite{Farb:mcg}) that the entire orbifold fundamental group of ${\cal
M}_{g,n}$ is carried by its end; in particular the fundamental group
of the end is not virtually nilpotent.  This argument actually shows
that no finite cover of ${\cal M}_{g,n}$ admits a complete, finite
volume Riemannian metric of pinched negative curvature.

\section{Quasi-flats in the pants complex}
\label{section:rank}

Let $S$ be a surface of genus $g$ with $n$ boundary components.  
A {\em subsurface} $R \subset S$ is a compact connected embedded
surface lying in $S$.  The subsurface $R$ is {\em essential} if its
boundary components are homotopically essential in $S$.

We say $S$ {\em decomposes} into essential subsurfaces
$R_1, \ldots, R_k$ if each $R_j$ may be modified by an isotopy so that
they are pairwise disjoint and $S - R_1 \disjunion \ldots \disjunion
R_k$ is a collection of open annular neighborhoods of simple closed
curves on $S$, each isotopic to a boundary component of $R_j$.

Let $r(S)$ denote the maximum number $k$ in any decomposition of $S$
into essential subsurfaces $R_1,
\ldots, R_k$ such that each $R_j$,
$j=1,\ldots, k$ has genus at least one, or at least four boundary
components.  Then $r(S)$ is greatest integer less which is at most
$(d(S)+1)/2.$

\begin{theorem}
The complex $C_{\bf P}(S)$ contains a quasi-flat of dimension $r(S)$.
\end{theorem}

\bold{Proof:}  
The surface $S$ decomposes into subsurfaces
$$R_1, \ldots, R_{r(S)}, T$$ so that $d(R_j) = 1$ for each
$j$ and either $T$ is empty or $d(T) = 0$.  
Let $n = r(S)$.  We describe a quasi-isometric
embedding of the Cayley graph for $\zed^n$ with the standard
generators into the complex $C_{\bf P}(S)$.  

Let $c_j$ be a vertex in the curve complex $\calC(R_j)$.
Then together with the core curves of the open annuli
in $S - R_1 \disjunion \ldots \disjunion R_{r(S)} \disjunion T$, 
the curves $c_j$ form a pants decomposition $P = P(c_1,
\ldots, c_n)$ of $S$.

We let $g_j \colon \zed \to \calC(R_j)$ be a geodesic so that
$g_j(0) = c_j$.  Note that $\calC(R_j)$ is the Farey graph, so we may
take any bi-infinite geodesic $g_j$ in $\calC(R_j)$.

We claim that the embedding
$$Q \colon \zed^n \to C_{\bf P}(S)$$
defined by 
$$Q(k_1, \ldots, k_n)  = P(g_1(k_1), \ldots, g_n(k_n))$$
is a quasi-isometry, whose constants depend only on $S$.

Let $\vec{k} = (k_1,\ldots, k_n)$ and $\vec{l} = (l_1, \ldots, l_n)$.
Since elementary moves along $g_j$ can be made independently in each
$R_j$, we have
$$d_{\bf P}(Q(\vec{k}), Q(\vec{l})) \le \sum_{j=1}^n 
|l_j - k_j|=d_{\zed^n}(\vec{k},\vec{l})$$
which shows that $Q$ is $1$-Lipschitz.

Given $R_j$, the 
projection $\pi_{R_j}(Q(\vec{k}))$ to $R_j$ described in
\S\ref{subsection:cc} simply
picks out the curve $g_j(k_j)$ so we have
$$\pi_{R_j}(Q(\vec{k})) = g_j(k_j).$$ 
Thus, the projection distance
$$d_{R_j}(Q(\vec{k}), Q(\vec{l})) = 
d_{R_j}(g_j(k_j), g_j(l_j))$$
which is simply $|k_j - l_j|$ since $g_j$ is a geodesic in
$\calC(R_j)$.

By Theorem 6.12 of \cite{Masur:Minsky:CCII}, there exists $M_0=M_0(S)$ so that for all
$M\geq M_0$ there exist constants $K_0$ and $K_1$ 
so that if we let $P_{\vec{k}} = Q(\vec{k})$ and $P_{\vec{l}} =
Q(\vec{l})$ then we have the inequality $$ \sum_{\stackrel{Y \subseteq
S}{d_Y(\pi_Y(P_{\vec{k}}), \pi_Y(P_{\vec{l}})) >M}} 
d_Y(P_{\vec{k}},P_{\vec{l}})
\le K_0 d_{\bf P}(P_{\vec{k}}, P_{\vec{l}}) + K_1.$$
But the left-hand-side of the inequality is bounded below by
$$\max_j |k_j - l_j| \ge \frac{\sum_j |k_j - l_j|}{n}.$$
Thus, $Q$ is a quasi-isometric embedding.
\qed

\section{The case of low genus}
\label{section:hierarchies}

In the case where $d(S) = 1$, where $\interior(S)$ is homeomorphic to
a punctured torus or four-times-punctured sphere, the pants complex
$C_{\bf P}(S)$ is identified with the curve complex $\calC(S)$, which
is Gromov-hyperbolic (see \cite{Masur:Minsky:CCI}). 

Together with the previous section, this observation leaves one case
unattended, namely that when $d(S) = 2$.  In this case, $\interior(S)$
is homeomorphic either to a doubly-punctured torus, or a
five-times-punctured sphere.  In this section we prove the following.

\begin{theorem}
Let $S$ be such that $d(S)$ equals $1$ or $2$.  Then the
Weil-Petersson metric  on $\Teich(S)$ is Gromov-hyperbolic.
\label{theorem:doubly}
\end{theorem}

The case $d(S) = 1$ is proven in \cite{Masur:Minsky:CCI}, since
$\calC(S) = C_{\bf P}(S)$ in this case.  Thus we are left to treat the
case when $d(S) = 2$.
In this case, we apply results of the second author which, although
initially phrased in the context of groups with their word metrics,
apply to general metric spaces.  

Our argument will employ the notion of {\em relative hyperbolicity}
developed in \cite{Farb:bcp}.  In essence, a metric space is {\em relatively
hyperbolic} relative to a collection of subsets if the result of
crushing those subsets to have diameter 1 is a hyperbolic metric space.
In the case $d(S) = 2$, we will show that $C_{\bf P}(S)$ is relatively
hyperbolic relative to regions consisting of pants decompositions
containing a single curve.  Since, in this case, such regions are
themselves hyperbolic, it is possible to establish
Gromov-hyperbolicity $C_{\bf P}(S)$ by showing that paths in $C_{\bf
P}(S)$ that determine quasi-geodesics in the relative space satisfy certain boundedness properties with respect
to their trajectories through these regions (this is the {\em bounded
region penetration} property, below).

\bold{Remark:}  One can prove theorem~\ref{theorem:doubly} by directly
demonstrating a thin-triangles condition after replacing geodesics in
$C_{\bf P}(S)$ by the {\em hierarchies} of \cite{Masur:Minsky:CCII}
(this was our original approach to the argument).  We have chosen
instead to employ the theory relative hyperbolicity as it is 
more familiar, and unifies these cases with the higher genus cases.  
Indeed, when $d(S) >2$ the natural regions with respect to which
$C_{\bf P}(S)$ is relatively hyperbolic (sub-graphs of pairs of pants
containing a given curve $\alpha$) are {\em not} themselves hyperbolic;
they are the quasi-flats of the previous section.  
\medskip

\bold{Relative hyperbolicity.}  Let $(X,d)$ be a geodesic metric
space, and let $H_\alpha$ be a 
collection of connected subsets of $X$, with index $\alpha$ in an
index set $A$. Then the {\em electric space} $\widehat{X}$ relative to
the regions $\{H_\alpha\}$ is obtained by collapsing each region to
have diameter one, as follows (see \cite[Sec. 3]{Farb:bcp}).  Adjoin
to the space $X$ a single point $c_\alpha$ for each $\alpha \in A$ by
connecting $c_\alpha$ to each point of $H_\alpha$ by a segment of
length $1/2$.  Let $\widehat{X}$ denote the resulting path-metric
space and let $d_e(.,.)$ denote path distance in $\widehat{X}$.  Given
a path $w$ in $X$ we obtain a path in $\widehat{X}$ by replacing
segments where $w$ travels in $H_\alpha$ with a path joining the
endpoints of the segment to $c_\alpha$.  As in \cite{Farb:bcp}, we
denote this path-replacement procedure by $X \to \widehat{X}$ or $w
\mapsto \hat{w}$.  We denote by $I(w)$ the {\em initial point} of $w$
and by $T(w)$ the {\em terminal point} of $w$.  The points $I(w)$ and
$T(w)$ depend on the choice of parameterization.

If $\hat{w}$ is a ($k$-quasi) geodesic in $\widehat{X}$ we say $w$ is
a {\em relative ($k$-quasi) geodesic } in $X$.  If a path $w$ in $X$
(or $\hat{w}$ in $\widehat{X}$) passes through some region $H_\alpha$
we say it {\em penetrates} $H_\alpha$.  A path $w \in X$ (or $\hat{w}$
in $\widehat{X}$) has {\em no backtracking} if for every region
$H_\alpha$ that $\hat{w}$ penetrates, once it leaves $H_\alpha$ it
never returns.  The space $X$ is {\em hyperbolic relative to
$\{H_\alpha\}_{\alpha \in A}$} if the electric space $\widehat{X}$ is
Gromov hyperbolic.

\medskip

In the pants
complex, consider following collection of regions: for each $\alpha
\in \calC(S)$ let
$$H_\alpha = \{ P \in
C_{\bf P}(S) \st \alpha \in P \}.$$
Then we have the following theorem (cf. \cite[Thms. 1.2,
1.3]{Masur:Minsky:CCI}).
\begin{lem}
The graph $C_{\bf P}(S)$ is hyperbolic relative to the regions $\{H_\alpha\}$.
\end{lem}
\bold{Proof:}  It suffices to show that the electric space
$\widehat{C_{\bf P}(S)}$ with 
respect to the regions $\{H_\alpha \}$ 
is quasi-isometric to the curve
complex $\calC(S)$, which is 
Gromov-hyperbolic by \cite[Thm. 1.1]{Masur:Minsky:CCI}. 

To see this, let $\Gamma = C_{\bf P}(S)$, let 
$\widehat{\Gamma}$ be the electric space associated to the regions
$\{H_\alpha\}$, and let $c_\alpha$ be the point added to $\Gamma$ at
distance 1/2 from each point of $H_\alpha$ to form $\widehat{\Gamma}$.
Consider the mapping $$q \colon \calC^0(S) \to \widehat{\Gamma}$$ from
the zero-skeleton of $\calC(S)$ to $\Gamma$ obtained by setting
$q(\alpha) = c_\alpha$.  Note that given a pants decomposition $P$, if
$\beta$ is an element of $P$ then $P$ lies a distance $1/2$ from
$c_\beta$, so the image $q(\calC^0(S))$ is $1/2$-dense.

Moreover, we have
$d_{\calC(S)}(\alpha,\beta) = 1$ if and only if there is a $P$ for
which $\alpha \in P$ and $\beta \in P$.  But $\alpha \cup \beta
\subset P$ holds if and only if
the regions $H_\alpha$ and  $H_\beta$ intersect, which holds if and
only if
$$d_{\widehat{\Gamma}}(c_\alpha,c_\beta)= 1.$$ 
Thus, the map $q$ is 1-bi-Lipschitz, and since the image is $1/2$-dense
we may construct a 2-Lipschitz inverse to $q$.  Thus, $q$ a quasi-isometry.
\qed

\bold{Bounded region penetration.}  We now recall results of
\cite{Farb:bcp} detailing a criterion on relative quasi-geodesics
that will serve to ensure hyperbolicity of $C_{\bf P}(S)$ when $d(S)
=2$.  The following definition is analogous to the ``bounded coset
penetration property'' in \S 3.3 of \cite{Farb:bcp}.

\begin{defn}[bounded region penetration]
The pair $(X,\{H_\alpha\})$ satisfies the {\em bounded region
penetration} property if, for every $P \ge 1$ there is a constant $c =
c(P) >0$ so that if $u$ and $w$ are relative $P$-quasi-geodesics
without backtracking so that the initial and terminal points of $u$
and $w$ satisfy $d_X(I(u),I(w))\le 1$ and $d_X(T(u),T(w)) \le 1$ then
the following holds:
\begin{enumerate}
\item If $u$ penetrates a region $H_\alpha$ but $w$ does not penetrate
$H_\alpha$, then $u$ travels an $X$-distance at most $c$ in
$H_\alpha$.
\item If both $u$ and $w$ penetrate a region $H_\alpha$ then the
points at which $u$ and $w$ first enter $H_\alpha$ lie an
$X$-distance at most $c$ from one another and likewise for the exit points.
\end{enumerate}
\end{defn}
An important application is the
following theorem in which bounded region penetration is used to
bootstrap from  hyperbolicity relative to hyperbolic regions to
hyperbolicity of the original metric space.  

\begin{theorem}
Suppose $X$ is hyperbolic relative to the regions $\{H_\alpha\}$ and
that the pair $(X,\{H_\alpha\})$ has the bounded region penetration
property.  Then if the regions 
$H_\alpha$ are themselves $\delta$-hyperbolic metric spaces for some
$\delta >0$, then $X$ is a Gromov-hyperbolic metric space.
\label{theorem:hyp:rel:hyp}
\end{theorem}

\bold{Proof:} {\em (of Theorem \ref{theorem:hyp:rel:hyp}).} 
Theorem \ref{theorem:hyp:rel:hyp} is simply a recasting of the remark
following \cite[Thm. 3.8]{Farb:bcp} from groups to general metric
spaces. The proof works verbatim in this case, with the following
addition: one replaces the use of the theorem that linear
isoperimetric inequality for a group implies that the group is
Gromov-hyperbolic by the corresponding theorem for metric spaces with
a well-defined notion of area.  Such a theorem is proven by Bowditch
in \cite{Bowditch:ip}; in this case one can use for area the
combinatorial area 
of the simplicial complex whose $1$-skeleton is the pants graph
$C_{\bf P}(S)$ and whose $2$-cells consist of five types of loops with
$3$, $4$, $5$, 
and $6$ edges, and no other edges between vertices (this is a variant
of the $2$-complex studied by Hatcher-Thurston in
\cite{Hatcher:Thurston:pants}.  It is shown to be simply connected in
\cite[Thm. D]{Hatcher:pants}).  One may verify that Bowditch's proof
extends to the locally infinite case. 
\qed

\bold{Proof:}  {\em (of Theorem~\ref{theorem:doubly}).}  
The condition
$d(S) = 2$ implies that each pants decomposition of $S$ is built from
exactly two disjoint simple closed curves on $S$. 
\begin{lem}
Suppose $d(S) = 2$.  Then there is a $\delta$ so that for each $\alpha
\in \calC(S)$, the region $H_\alpha$ is $\delta$-hyperbolic.
\end{lem}

\bold{Proof:}  Given $\alpha \in \calC(S)$, let $Y_\alpha$ denote the
connected component of the complement of an embedded open annular
neighborhood of $\alpha$ for which $d(Y_\alpha) =1$.   Then the region
$H_\alpha$ is isometric to the curve complex
$\calC(Y_\alpha)$.  Again, \cite[Thm. 1.1]{Masur:Minsky:CCI}
implies that $H_\alpha$ is hyperbolic.
\qed

To prove Theorem~\ref{theorem:doubly}, then, it suffices to prove that
when $d(S) = 2$, the pair $(C_{\bf P}(S), \{H_\alpha\})$ has the
bounded region penetration property.  To this end, let $u$ and $w$ be
two relative $P$-quasi-geodesics in $C_{\bf P}(S)$ without
backtracking, so that $d_{\bf P}(I(u),I(w)) \le 1$ and $d_{\bf
P}(T(u),T(w)) \le 1$.  Let $H_\alpha$ be a region which $u$ penetrates
but $w$ does not. Being relative $P$-quasi-geodesics in the relatively
hyperbolic space $(C_{\bf P}(S),\{H_\alpha\})$, it follows form the
definition that the projections $\hat{u}$ and $\hat{w}$
$D$-fellow-travel in the electric space $\widehat{C_{\bf P}(S)}$ for
some $D>0$ depending only on $P$.  For simplicity of notation let
$\Gamma = C_{\bf P}(S)$ and let $\widehat{\Gamma}$ be the associated
electric space relative to the regions $\{ H_\alpha \}$.

We observe that in our circumstances, the path replacement $u \mapsto
\hat{u}$ 
can be viewed as producing an explicit path in $\calC(S)$ from a path
$u$ in $\Gamma$.  Since $d(S) =2$, each pants decomposition in of $S$
has two elements, so a path $u$ in $\Gamma$ is a sequence of edges in
$\calC(S)$ each of which is joined to the previous one at one of its
two endpoints.  Let $\tilde{u}$ denote the path in $\calC(S)$ obtained
by removing all but the first and last edges in any sequence of
consecutive edges in $u$ that all contain a single vertex.  The path
$\tilde{u}$ in $\calC(S)$ has image $\hat{u}$ under the quasi-isometry
$q$ (up to segments of length $1/2$ at the endpoints).

To employ the extra information the curve complex provides, we work
with $\tilde{u}$ and $\tilde{w}$ rather than $\hat{u}$ and $\hat{w}$.
The condition that $u$ is a relative $P$-quasi-geodesic without
backtracking simply means that the path $\tilde{u}$ is a
$P$-quasi-geodesic in $\calC(S)$ that never repeats a vertex.  
Proving that $(\Gamma,\{H_\alpha\})$ satisfies bounded region
penetration property, then, reduces to verifying that for paths $u$ and $w$
in $\Gamma$ whose for which $d_{\Gamma}(I(u),I(w)) \le 1$ and
$d_{\Gamma}(T(u),T(w)) \le 1$, and whose corresponding paths
$\tilde{u}$ and $\tilde{w}$ are $P$-quasi-geodesics without
backtracking we have:
\begin{enumerate}
\item[$1'$] If $\tilde{u}$ encounters a vertex $v_\alpha$ that $\tilde{w}$ avoids,
the vertices $v_\beta$ and $v_\gamma$ adjacent to $v_\alpha$ on
$\tilde{u}$ have distance
$$d_{Y_\alpha}(v_\beta,v_\gamma) < c$$ in the subsurface $Y_\alpha$.
\item[$2'$] If $\tilde{u}$ and $\tilde{w}$ each encounter a vertex
$v_\alpha$, then the vertices $v_\beta \in \tilde{u}$ and $v_{\beta'}
\in \tilde{w}$ just prior to
the encounter with $v_\alpha$ satisfy
$$d_{Y_\alpha}(v_\beta,v_{\beta'}) < c$$ and likewise for the points
$v_\gamma$ and $v_{\gamma'}$ on $\tilde{u}$ and $\tilde{w}$ just
following the encounter with $v_\alpha$.
\end{enumerate}
(For the remainder of this section we will denote by
$v_\alpha$ the vertex 
in $\calC(S)$ corresponding to the simple closed curve $\alpha$ on $S$
to avoid notational confusion).
To see property $(1')$ implies property $(1)$ above, note
that the condition that $\{v_\beta , v_\alpha, 
v_\gamma \}$ is a sub-segment of $\tilde{u}$ implies that
$P = \{v_\alpha,v_\beta\}$ and $P' = \{ v_\alpha, v_\gamma\}$ are the pants
decompositions along the path $u$ where $u$ enters and exits the
region $H_\alpha$.  If $d_{Y_\alpha}(v_\beta,v_\gamma) <c$ then there
is a sequence of pants decompositions joining $P$ to $P'$ of length at
most $c$ given by taking a geodesic 
$$\{v_\beta = v_0, \ldots, v_N = v_\gamma \}$$
joining $v_\beta$ to $v_\gamma$ in
$\calC(Y_\alpha)$ and taking the sequence of pants decompositions to be
$\{P_j = \{v_\alpha,v_j\}\}_j$.  One argues similarly that property
$(2')$ implies property $(2)$.

We now verify that properties $(1')$ and $(2')$ hold.
Since $\tilde{u}$ does not backtrack, we may choose points $x$ and $y$
on $\tilde{u}$ on either side of the vertex
$v_\alpha$ as follows.
Either $x$ is an endpoint of 
$\tilde{u}$ or $x$ lies at distance $2D$ in $\calC(S)$ from
the vertex $v_\beta$ adjacent on $\tilde{u}$ to $v_\alpha$ whichever is closer along $\tilde{u}$.
The point $y$ is either an endpoint of $\tilde{u}$ or $y$ lies at
distance $2D$ along $\tilde{u}$ from the vertex 
$v_\gamma$ adjacent to $v_\alpha$ whichever is closer along
$\tilde{u}$.  There is then a path $p$ on 
$\tilde{w}$ (which does not encounter $v_\alpha$) from the nearest
point to $x$ on $\tilde{w}$ to the nearest 
point to $y$.  Letting $p_x$ be the shortest path joining $x$ to
$\tilde{w}$ and letting $p_y$ be the shortest path joining $y$ to
$\tilde{w}$ in $\calC(S)$, the concatenation
$$q = p_x \compos p \compos p_y$$
is a path in $\calC(S)$ that avoids the vertex $v_\alpha$.  
Moreover, the path $q$ has length at most $8D$.  

Letting $q_x$ be the path along $\tilde{u}$ joining $x$ to 
$v_\beta$, and letting $q_y$ be the path along
$\tilde{u}$ joining $y$ to $v_\gamma$, we
have a path $$r = q_x \compos p \compos q_y$$ of length at most
$12D$ that avoids the vertex $v_\alpha$ entirely.    

The path $r$ describes a path in the curve complex $\calC(S)$ so that
each vertex along the interior of $r$ corresponds to a curve that
either lies in $\calC(Y_\alpha)$ or 
intersects $\bdry Y_\alpha$.  It follows from
lemma~\ref{lemma:Lipschitz:simplex} that
any two consecutive vertices $z$ and $z'$ on $r$ satisfy
$$d_{Y_\alpha}(z,z') \le 2.$$ 
Thus we have the bound $$d_{Y_\alpha}(v_\alpha,v_\beta) \le 24D.$$

\medskip

A similar argument proves property (2) of the bounded region
penetration property holds.  Choose a point $x$ on $\tilde{u}$ so that
either $x$ is the first vertex of $\tilde{u}$ or $x$ is at distance
$2D$ along $\tilde{u}$ from $v_\beta$, whichever is closer along
$\tilde{u}$ to $v_\beta$.  By fellow traveling, there is a point $y$
on $\tilde{w}$ and a path $p$ of length at most $D$ joining $x$ to $y$
in $\calC(S)$.  The path along $\tilde{w}$ joining $y$ to the last
vertex $v_\eta$ prior to $v_\alpha$ along $\tilde{w}$ has length at
most $2D(1 + P)$ in $\calC(S)$, so there is a path $r$ in $\calC(S)$
of total length bounded by $3D + 2D(1 + P)$ joining $v_\beta$ to
$v_\eta$ that does not hit $v_\alpha$ in its interior.  Again, by
lemma~\ref{lemma:Lipschitz:simplex} any two consecutive vertices $z$
and $z'$ on $r$ have the property that $$d_{Y_\alpha}(z,z') \le 2,$$
so we have $$d_{Y_\alpha}(v_\alpha, v_\eta) \le 2(3D + 2D(1+P)).$$ The
same argument proves that points on $v_\alpha'$ and $v_\eta'$ adjacent
to $v_\alpha$ where $\tilde{u}$ and $\tilde{w}$ depart from $v_\alpha$
also have bounded distance $d_{Y_\alpha}(v_\alpha',v_\eta').$

Having verified that properties $(1')$ and $(2')$ hold, we conclude
that the pair $(\Gamma, \{ H_\alpha \})$ has the bounded region 
penetration property.  The theorem follows from
theorem~\ref{theorem:hyp:rel:hyp}. 
\qed

\section{Questions}

We close the paper with some natural questions.
\begin{quest}
{\bf (McMullen)} Does ${\cal M}_{g,n}$ admit a complete, nonpositively
curved Riemannian metric?
\end{quest}
McMullen's K\"ahler metric on $\calM_{g,n}$ is complete but not
nonpositively curved, while the Weil-Petersson metric is nonpositively
curved but is not complete.  Is there a possible compromise?

\begin{quest}
What is the geometric rank of 
\begin{itemize}
\item the Weil-Petersson metric?
\item the Teichm\"uller metric?
\item the mapping class group?
\end{itemize}
\end{quest}
Theorem~\ref{theorem:rank} gives the lower bound $d(S)/2$ to the rank
of the Weil-Petersson metric, while \cite{Minsky:projections}
and \cite{Farb:Lubotzky:Minsky} establish the lower bound $d(S)$ for the rank of the 
Teichm\"uller metric and the mapping class group respectively.

The answers to these rank questions seem to be essential to understanding
quasi-isometric rigidity questions in Teichm\"uller space and the
mapping class group.

\begin{quest}
If $\interior(S)$ is homeomorphic to a doubly-punctured torus or
5-times-punctured sphere, are the sectional curvatures of the
Weil-Petersson metric bounded away from zero?
\end{quest}
Were the (geodesically convex) Weil-Petersson metric to have curvature
pinched from above by a negative constant, its Gromov-hyperbolicity in
this case would be an immediate consequence.

\noindent
{\sc \scriptsize Math Department,
University of Chicago,
5734 S. University Ave.,
Chicago, IL 60637\\
email: \ \  brock@math.uchicago.edu, \ \ farb@math.uchicago.edu}

\end{document}